\documentclass[11pt,onecolumn,draftcls]{IEEEtran}

\usepackage{verbatim}
\usepackage{amsfonts,amsmath,mathrsfs,amssymb,amsbsy}
\usepackage[final]{graphicx}
\usepackage{times,cite}
\usepackage{subfig}

\long\def\symbolfootnote[#1]#2{\begingroup%
\def\thefootnote{\fnsymbol{footnote}}\footnote[#1]{#2}\endgroup}

\usepackage{verbatim}
\usepackage[]{float,latexsym}
\usepackage{amsfonts,amsmath,mathrsfs,amssymb,amsbsy}
\usepackage{url}

\newtheorem{theorem}{Theorem}[section]
\newtheorem{corrly}{Corollary}[theorem]

\newtheorem{algorithm}{Algorithm}[section]
\newtheorem{assumption}{Assumption}[section]
\newtheorem{problem}{Problem}[section]

\newcommand{\Prob}{\mathbb{P}}
\newcommand{\Expect}{\mathbb{E}}

\long\def\symbolfootnote[#1]#2{\begingroup%
\def\thefootnote{\fnsymbol{footnote}}\footnote[#1]{#2}\endgroup}

\DeclareMathOperator*{\esssup}{ess\,sup}

\newcommand{\CADD}{{\mathsf{CADD}}}

\newcommand{\WADD}{{\mathsf{WADD}}}
\newcommand{\FAR}{{\mathsf{FAR}}}

\newcommand{\PDC}{{\mathsf{PDC}}}

\newcommand{\tauc}{\tau_{\scriptscriptstyle \mathrm{C}}}

\newcommand{\taugc}{\tau_{\scriptscriptstyle \mathrm{GC}}}
\newcommand{\taugd}{\tau_{\scriptscriptstyle \mathrm{GD}}}

\newcommand{\tauw}{\tau_{\scriptscriptstyle \mathrm{W}}}

\begin{document}
\title{Data-Efficient Minimax Quickest Change Detection with Composite Post-Change Distribution}
\author{Taposh Banerjee, \textit{Student Member, IEEE},  and Venugopal V. Veeravalli, \textit{Fellow, IEEE \vspace{-3ex}}
}
%\author{\IEEEauthorblockN{Taposh Banerjee and Venugopal V. Veeravalli}\\
%\IEEEauthorblockA{ECE Department and
%Coordinated Science Laboratory\\
%University of Illinois at Urbana-Champaign, Urbana, IL\\
%Email: banerje5,vvv@illinois.edu}}

%\onecolumn
\maketitle

%\ninept

\symbolfootnote[0]{\small
This research was supported by the National Science Foundation (NSF) under grant DMS 12-22498 and by the Defense Threat Reduction Agency (DTRA)
under subcontract 147755 at the University of Illinois, Urbana-Champaign from prime award HDTRA1-10-1-0086.
Preliminary version of this paper has been presented at 2014 IEEE International Symposium
on Information Theory.

The authors are with the Department of Electrical and Computer Engineering and Coordinated Science Laboratory,
University of Illinois Urbana Champaign, Urbana, IL 61801 USA
(e-mail: banerje5@illinois.edu; vvv@illinois.edu)
}

\vspace{-0.3cm}

\begin{abstract}
The problem of quickest change detection is studied, where there is an additional
constraint on the cost of observations used before the change point and where the
post-change distribution is composite. Minimax formulations are proposed for this problem.
It is assumed that the post-change family of distributions has a member which is
least favourable in some sense.
An algorithm is proposed in which on-off observation control is employed using the least
favourable distribution, and a generalized likelihood ratio based approach is used for change detection.
Under the additional condition
that either the post-change family of distributions is finite, or both the pre- and post-change distributions
belong to a one parameter exponential family,
it is shown that the proposed algorithm is asymptotically optimal, uniformly for all possible
post-change distributions.
\end{abstract}
%\tb{\textbf{Need to argue that the MCuSum algorithm is equivalent to a GLRT}}

\begin{keywords}
Asymptotic optimality, CuSum, exponential family, generalized likelihood ratio, least favourable distribution, minimax, observation control, quickest change detection,  unknown post-change distribution.
\end{keywords}

%\vspace{-0.1cm}
\section{Introduction} \label{sec:Intro}
%\vspace{-0.1cm}
The problem of detecting an abrupt change in the statistical properties of a measurement
process is encountered in many engineering applications. Applications include detection
of the appearance of a sudden fault/stress in a system being monitored, e.g., bridges,
historical monuments, power grids, bird/animal habitats, etc. Often in these
applications the decision making has to be done in real time, by taking measurements
sequentially. In statistics this detection problem is formulated within the framework of quickest
change detection (QCD) \cite{veer-bane-elsevierbook-2013}, \cite{bane-chen-garcia-veer-icassp-2014}.

In the QCD problem, the objective is to detect an abrupt change in the distribution of a sequence
of random variables. The random variables follow a particular distribution in the beginning, and
after an unknown point of time, follow another distribution.
This problem is well studied in the literature \cite{veer-bane-elsevierbook-2013}, \cite{poor-hadj-qcd-book-2009}, \cite{tart-niki-bass-2014}.
%In the various popular problem formulations studied in the literature,
The objective is to find a stopping time for the random variables so as to minimize a suitable
metric on the average detection delay subject to a constraint on a suitable metric on the false alarm rate.
When the pre- and post-change distributions are known, the optimal stopping rule,
for all the popular QCD formulations in the literature, is a single threshold test, where a sequence of statistics
is computed using the likelihood ratio of the observations, and a change is declared the first time
the sequence of statistics crosses a threshold. The threshold is chosen to meet the constraint on
the false alarm rate.
For example, a popular algorithm in the literature that has some strong optimality properties is the Cumulative Sum (CuSum) algorithm
 (see Section~\ref{sec:BkGrdQCDDEQCD} for a precise statement).   
In the CuSum algorithm, the cumulative log likelihood ratio of the observations is computed over time. If the 
accumulated statistic is below zero, it is reset to zero. A change is declared when the 
accumulated statistic is above a threshold. 
 
In practice, often the post-change distribution is not known or known only to belong to a
parametric family of distributions.
Moreover, the change occurs rarely and
it is of interest to constrain the number of observations (data) used before the change point.

The classical problem of detecting a change when the post-change distribution is unknown
(and with no observation control) has been well studied
in the literature. In the parametric setting, where the post-change distribution is assumed to belong to
a parametric family, there are three main approaches: generalized likelihood ratio (GLR) based, mixture based
and adaptive estimates based approaches.
In the nonparametric setting, one approach has been to take a robust approach to the QCD problem.
See \cite{veer-bane-elsevierbook-2013}, \cite{poor-hadj-qcd-book-2009}
and \cite{tart-niki-bass-2014} for a review.

In \cite{bane-veer-sqa-2012} and \cite{bane-veer-IT-2013} we studied the classical QCD problems with
an additional constraint on a suitable metric for the cost of observations used before the change point.
We called these formulations data-efficient quickest change detection (DE-QCD).
For the case when the pre- and post-change distributions are known, 
we showed that two-threshold generalizations of the classical single-threshold QCD tests are asymptotically
optimal for the proposed formulations.
Specifically, in the two-threshold test, there are two thresholds. A sequence of statistics is computed over time
using the likelihood ratio of the observations. A change is declared the first time
the sequence of statistics crosses the larger of the two thresholds. If the computed statistic is below the
upper threshold, then the next observation is taken only if the statistic is above the smaller
of the two thresholds. The upper threshold is used to control the false alarm rate,
and the lower threshold is chosen to control the cost of observations before the change point.
For example, we proposed an algorithm called
the data-efficient cumulative sum (DECuSum) algorithm in \cite{bane-veer-IT-2013}, which is a two-threshold 
generalization of the CuSum algorithm. In the DECuSum algorithm also
the cumulative log likelihood ratio of the observations is computed over time. However, when the accumulated statistic
goes below zero, instead of resetting it to zero, the undershoot of the statistic is exploited for skipping consecutive samples.
Thus, the likelihood ratio of the observations is used for data-efficiency as well as for stopping.  
However, if the post-change distribution is not known, then it is not clear 
what statistic should be used for skipping samples for data-efficiency, while 
at the same time detecting the change in an optimal manner. 

In this paper we combine the ideas from \cite{bane-veer-IT-2013} and from the QCD literature
for the case where the post-change distribution is unknown to study DE-QCD problems when the post-change distribution is composite.
We assume that the post-change family of distributions has
a \textit{least favorable} member (see Assumption~\ref{assum:leastfav} for a precise definition). 
Based on this assumption, we propose an algorithm
called the generalized data-efficient cumulative sum (GDECuSum) algorithm.
In this algorithm on-off observation control is performed using the DECuSum algorithm
designed for the least favorable distribution,
and the change is detected using a GLRT based CuSum algorithm; the latter is called the GCuSum algorithm in the following. 
The GCuSum algorithm, studied in \cite{lord-amstat-1971} and \cite{tart_polu_ISIF_2008}, is a GLRT based 
extension of the CuSum algorithm from \cite{page-biometrica-1954}.  
Thus, the GDECuSum algorithm is an extension of the GCuSum algorithm with the feature of on-off observation 
control introduced to control the cost of observations used before the change point. 

We provide a detailed performance analysis of the GDECuSum algorithm. 
The performance analysis reveals (see Section~\ref{sec:GDECuSum_Perf} for mathematically precise statements)
that the false alarm rate of the 
GDECuSum algorithm is as good as the false alarm rate of the GCuSum algorithm. 
Also, the delay of the GDECuSum algorithm is within a constant of the delay of the
GCuSum algorithm. We will show that these two results on the delay and false alarm analysis 
can be used to prove the asymptotic optimality of the GDECuSum algorithm for the proposed formulations, 
whenever the GCuSum algorithm is asymptotically optimal for the classical QCD formulations. 
The GCuSum algorithm is asymptotically optimal for the classical formulations, for example, 
for the following cases: (i) when the post-change family of distributions is finite, and (ii) if both
the pre- and post-change distributions belong to a one-parameter
exponential family. 

The assumption that the post-change distribution
belongs to a finite set of distributions is satisfied in many practical applications.
For example, it is satisfied in
the problem of detecting a power line outage in a power grid \cite{bane-chen-garcia-veer-icassp-2014},
or in a multi-channel scenario where the observations are vector valued and a change affects the distribution of only
a subset of the components (each component for example may correspond
to the output of a distinct sensor on a sensor board) \cite{tart-veer-fusion-2002}, \cite{mei-biometrica-2010}. Also,
see \cite{tart_polu_ISIF_2008} for a possible scenario.

The paper is organized as follows.
In Section~\ref{sec:ProbForm} we propose a modified version of the minimax problem formulations from \cite{bane-veer-IT-2013}.
In Section~\ref{sec:BkGrdQCDDEQCD} we provide a brief review of QCD and DE-QCD relevant to this paper.
In Section~\ref{sec:G-DE-CuSum} we propose the main algorithm of the paper, the GDECuSum algorithm.
In Section~\ref{sec:GDECuSum_Perf} we analyze the performance of the GDECuSum algorithm and 
discuss its optimality properties. 
 In Section~\ref{sec:Discussion} we discuss possible extensions of this work to 
 mixture based tests. We also discuss the case when a least favorable distribution does not exist.
 Finally, we discuss extensions of our work to window limited GLR tests from \cite{lai-ieeetit-1998}.  
In Section~\ref{sec:DECompositeNumerical} we compare the performance of the GDECuSum algorithm
with the approach of fractional sampling, in which the GCuSum algorithm is used to detect the change
and the constraint on the cost of observations is satisfied by skipping samples randomly, independent of the
observation process.
In Section~\ref{sec:Conclusions} we conclude the paper.

\section{Problem Formulation}
\label{sec:ProbForm}
A sequence of random variables $\{X_n\}$ is being observed. Initially, the random variables are i.i.d. with
p.d.f. $f_{0}$. At time $\gamma$, called the change point, the density of the 
random variables changes from $f_0$ to $f_{\theta}$, $\theta \in \Theta$. That is, we assume
that the post-change distribution belongs to a parametric family of distributions parameterized by $\theta$.
Both $\theta$ and $\gamma$ are unknown.
%We assume that
%\[\theta \in \Theta = \{\theta_1, \theta_2, \cdots, \theta_M\}.\]
We assume that $f_0 \neq f_\theta$ for all $\theta \in \Theta$.
%Thus, we assume that the unknown post-change distribution belongs to a finite set of distributions.
We denote by $\Prob^{\theta}_\gamma$ the underlying probability measure which governs such a sequence. We
use $\Expect^{\theta}_\gamma$ to denote the expectation with respect to this probability measure.
We use $\Prob_\infty$ ($\Expect_\infty$) to denote the probability measure (expectation) when the change
never occurs, i.e., the random variable $X_n$ has p.d.f. $f_0$, $\forall n$.
%We wish to detect this change in distribution as quickly as possible subject to a constraint on the false alarm rate.

In the classical QCD problem the objective 
is to detect the change in distribution as quickly as possible, subject to a constraint on the false alarm rate.
Since in the classical QCD there is no constraint on the cost of observations used before the change point,  
the optimal trade-off between delay and false alarm rate is achieved by utlizing 
all the observations for decision making. 

In many applications the change occurs rarely, corresponding to a large $\gamma$. As a result,
we also wish to control the number of observations used for decision making before $\gamma$.
We are interested in control policies involving causal three-fold decision making at each time step. Specifically,
based on the information available at time $n$, a decision has to be made whether to declare a change
or to continue taking observations. If the decision is to continue, then a decision has to be made whether
to use or skip the next observation for decision making

Mathematically, let $S_n$ be the indicator random variable defined as 
\begin{equation*}
\begin{split}
               S_{n} & =  \begin{cases}
               1 & \mbox{ if } X_n \mbox{ used for decision making }\\
               0 & \mbox{ otherwise.}
               \end{cases}
               \end{split}
\end{equation*}
The information available at time $n$ is denote by
\[\mathcal{I}_n = \{X_1^{(S_1)}, \cdots, X_n^{(S_n)}\},\]
where $X_k^{(S_k)}=X_k$ if $S_k=1$, else $X_k$ is absent from $\mathcal{I}_n$,
and
\[S_n = \phi_n (\mathcal{I}_{n-1}).\]
Here, $\phi_n$ denotes the control map.
Let $\tau$ be a stopping time for the sequence $\{\mathcal{I}_n\}$. A control policy is the collection
\[\Psi = \{\tau, \phi_1, \cdots, \phi_{\tau}\}.\]

We now propose two stochastic optimization problems where the objective is to minimize a metric
on delay, subject to constraints on a metric on the false alarm rate and a metric on the cost of
observations used before the change point $\gamma$.
We seek policies of type $\Psi$ to solve the proposed stochastic optimization problems.

We now define the metrics to be used in the problem formulations.
For delay we choose the following conditional average detection delay metric ($\CADD$) of Pollak \cite{poll-astat-1985}:
\begin{equation}
\label{def:CADD}
\CADD^{\theta}(\Psi) := \sup_{\gamma\geq 1} \; \Expect^{\theta}_\gamma[\tau-\gamma| \tau \geq \gamma].
\end{equation}
Note that the $\CADD$ is a function of the post-change parameter $\theta$.

For false alarm we choose the metric of false alarm rate ($\FAR$) used by Lorden in \cite{lord-amstat-1971} and by Pollak in \cite{poll-astat-1985}:
\begin{equation}
\label{def:FAR}
\FAR(\Psi) := \frac{1}{\Expect_\infty[\tau]}.
\end{equation}

To capture the cost of observations used before $\gamma$, we use the following variation of the
duty cycle metric proposed in \cite{bane-veer-IT-2013}, the Pre-change Duty Cycle ($\PDC$) metric:\footnote{The definition of $\PDC$ used in \cite{bane-veer-IT-2013} has an extra conditioning on $\{\tau \geq \gamma\}$.}
\begin{equation}
\label{def:PDC}
\begin{split}
 \PDC(\Psi)   &:= \limsup_{\gamma\to \infty} \; \Expect^{\theta}_\gamma\left[\frac{1}{\gamma}\sum_{n=1}^{\gamma-1} S_n \right] \\
 & = \limsup_{\gamma\to \infty} \; \Expect_\infty\left[\frac{1}{\gamma}\sum_{n=1}^{\gamma-1} S_n \right].
\end{split}
\end{equation}
%For simplicity we call this metric $\PDC$ as well.
Note that both the $\FAR$ and the $\PDC$ are \textit{not} a function of the post-change parameter $\theta$.

The first problem that we are interested in is the following:
\begin{problem}\label{prob:Pollak_Unknownf1}
\begin{eqnarray}
 \min_\Psi        &&   \CADD^{\theta}(\Psi) \nonumber\\
 \mbox{subj. to}  &&   \FAR(\Psi) \leq \alpha, \nonumber\\
\mbox{and}        &&  \PDC(\Psi)   \leq \beta,   \nonumber
\end{eqnarray}
where $0 \leq  \alpha, \beta \leq 1$ are given constraints.
\end{problem}

We are also interested in the problem where the $\CADD$ in Problem~\ref{prob:Pollak_Unknownf1} is replaced by the following
worst case average detection delay ($\WADD$) metric of Lorden \cite{lord-amstat-1971},
\begin{equation}
\label{def:WADD}
\WADD^{\theta}(\Psi) := \sup_{\gamma\geq 1} \; \esssup \Expect^{\theta}_\gamma[(\tau-\gamma)^+| \mathcal{I}_{\gamma-1}],
\end{equation}
where $x^+ := \max\{0, x\}$:
\begin{problem}\label{prob:Lorden_Unknownf1}
\begin{eqnarray}
 \min_\Psi        &&   \WADD^{\theta}(\Psi) \nonumber\\
 \mbox{subj. to}  &&   \FAR(\Psi) \leq \alpha, \nonumber\\
\mbox{and}        &&  \PDC(\Psi)   \leq \beta,   \nonumber
\end{eqnarray}
where $0 \leq  \alpha, \beta \leq 1$ are given constraints.
\end{problem}

For any policy $\Psi$ we have
\begin{equation}
\label{eq:CADDleqWADD}
\CADD^{\theta}(\Psi) \leq \WADD^{\theta}(\Psi).
\end{equation}

Our objective is to find an algorithm
that is a solution to both Problem~\ref{prob:Pollak_Unknownf1} and Problem~\ref{prob:Lorden_Unknownf1} 
uniformly for each $\theta \in \Theta$. However, it is not clear if such a solution
exists, even with $\beta=1$. As a result we seek a solution that is asymptotically optimal, for a given $\beta$, for each $\theta$, as $\alpha \to 0$.

In the rest of the paper we use $D(f_\theta \; ||\; f_0)$ and $D( f_0 \; ||\; f_\theta )$ to denote
\begin{eqnarray*}
D(f_\theta \; ||\; f_0) &:=& \Expect^{\theta}_1 \left[\log \frac{f_\theta(X_1)}{f_0(X_1)}\right],\\
D(f_0\; ||\; f_\theta ) &:=& - \Expect_\infty \left[\log \frac{f_\theta(X_1)}{f_0(X_1)}\right].
\end{eqnarray*}
We assume throughout that both $D(f_\theta \; ||\; f_0)$ and $D(f_0\; ||\; f_\theta )$ are finite and positive. 

\section{Classical QCD with Unknown Post-Change Distribution}
\label{sec:BkGrdQCDDEQCD}
In this section we review the results from \cite{lord-amstat-1971}, \cite{tart_polu_ISIF_2008} and \cite{bane-veer-IT-2013}
that are relevant to this paper.

We first review the lower bound on the performance of any test for an $\FAR$ of $\alpha$. 
Let
\[\Delta_\alpha := \{\Psi: \FAR(\Psi) \leq \alpha\}.\]
When the post-change density is $f_\theta$, a universal lower bound on the $\CADD^{\theta}$ over the class $\Delta_\alpha$ is given by (see \cite{lai-ieeetit-1998})
\begin{equation}
\label{eq:LB}
\inf_{\Psi \in \Delta_\alpha} \CADD^{\theta}(\Psi) \geq \frac{|\log \alpha|}{D(f_\theta \; ||\; f_0)}(1+o(1)) \mbox{ as } \alpha \to 0.
\end{equation}
By \eqref{eq:CADDleqWADD}, this is a lower bound on $\WADD^{\theta}$ as well.

\subsection{QCD with No Observation Control ($\beta=1$), $\theta$ Known}
We first consider the case when the post-change distribution is known to be $f_\theta$, i.e., when the post-change
parameter $\theta$ is known, and when there is no observation control, i.e., when $\beta=1$, in Problem~\ref{prob:Pollak_Unknownf1} and Problem~\ref{prob:Lorden_Unknownf1}.
Then the lower bound \eqref{eq:LB} is achieved by the cumulative sum (CuSum) algorithm \cite{page-biometrica-1954}, \cite{lord-amstat-1971}.
The CuSum algorithm is defined as follows:
\begin{equation}
\label{eq:CUSUM_algo}
\begin{split}
%C_0(\theta) &= 0,\\
C_{n}(\theta) &= \max_{1\leq k \leq n+1} \sum_{i=k}^n \log \frac{f_\theta(X_i)}{f_0(X_i)} \quad \mbox{ for } n\geq 1,\\
\tauc(\theta) &= \inf\{n \geq 1: C_n(\theta) \geq A\}.
\end{split}
\end{equation}
The statistic $C_n(\theta)$ can be computed recursively:
\begin{equation}
\label{eq:CUSUM_recursion}
\begin{split}
C_0(\theta) &= 0,\\
C_{n}(\theta) &= \left(C_{n-1}(\theta) +  \log \frac{f_\theta(X_n)}{f_0(X_n)}\right)^+ \mbox{ for } n\geq 1.
\end{split}
\end{equation}

The CuSum algorithm is asymptotically optimal for both Problem~\ref{prob:Pollak_Unknownf1} and Problem~\ref{prob:Lorden_Unknownf1} (with $\theta$ known and $\beta=1$)
due to \eqref{eq:CADDleqWADD} and because of the following result:
setting $A=\log 1/\alpha$ in \eqref{eq:CUSUM_algo} ensures that \cite{lord-amstat-1971}
\begin{equation}
\label{eq:CuSumPerf}
\begin{split}
\FAR(\tauc(\theta)) &\leq \alpha,\\
\WADD^{\theta}(\tauc(\theta)) &\leq \frac{|\log \alpha|}{D(f_\theta \; ||\; f_0)}(1+o(1)) \mbox{ as } \alpha \to 0.
\end{split}
\end{equation}
We note that the $\PDC$ of the CuSum algorithm is equal to 1.

\subsection{QCD with No Observation Control ($\beta=1$), $\theta$ Unknown}
We next consider the case when the post-change distribution is unknown, i.e., when the post-change
parameter $\theta$ is unknown, and again there is no observation control, i.e., $\beta=1$ in
Problem~\ref{prob:Pollak_Unknownf1} and Problem~\ref{prob:Lorden_Unknownf1}.
A natural extension of the CuSum algorithm for this case is the generalized likelihood ratio based CuSum algorithm.
We refer to the algorithm as the GCuSum algorithm and it is defined as follows:
\begin{equation}
\label{eq:GLRTCUSUM_algo}
\begin{split}
G_{n} &= \max_{1\leq k \leq n} \; \sup_{\theta \in \Theta'(\alpha)} \; \sum_{i=k}^n \log \frac{f_\theta(X_i)}{f_0(X_i)} \quad \mbox{ for } n\geq 1,\\
\taugc &= \inf\{n \geq 1: G_n \geq A\},
\end{split}
\end{equation}
where $\Theta'(\alpha)\subset \Theta$ can be a function of $\alpha$,
and is either equal to $\Theta$, or is allowed to be arbitrarily close and grow to $\Theta$ as $\alpha \to 0$.
The GCuSum algorithm has the following interpretation.
To detect a change when the post-change parameter is unknown,
a family of CuSum algorithms are executed in parallel,
one for each post-change parameter. A change is declared the first time a change
is detected in any one of the CuSum algorithms.
It can be shown that
\begin{equation}
\label{eq:CADDWADDWorstGamma_MCuSum}
\begin{split}
\WADD^{\theta}(\taugc) = \CADD^{\theta}(\taugc)=\Expect^{\theta}_1[\taugc-1].
\end{split}
\end{equation}
We also note that the $\PDC$ of the GCuSum algorithm is equal to 1.

The asymptotic optimality of the GCuSum algorithm is known for example in the following two cases:
when the post-change family is finite \cite{tart_polu_ISIF_2008},
and when the pre- and post-change distributions belong to
a one-parameter exponential family \cite{lord-amstat-1971}.

When the post-change set $\Theta$ is finite, i.e,
\[
\Theta = \{\theta_1, \cdots, \theta_M\},
\]
the GCuSum algorithm with $\Theta'(\alpha) = \Theta$ reduces to the following algorithm
\begin{equation}
\label{eq:GLRTCuSum_finitetheta_1}
\begin{split}
\taugc = \inf\left\{n \geq 1: \max_{1 \leq k \leq M} C_n(\theta_k) \geq A\right\},
\end{split}
\end{equation}
where $C_n(\theta_k)$ is the CuSum statistic \eqref{eq:CUSUM_recursion}
evaluated for $\theta = \theta_k$. Equation \eqref{eq:GLRTCuSum_finitetheta_1} can also be written as
\begin{equation}
\label{eq:GLRTCuSum_finitetheta_2}
\begin{split}
\taugc = \min_{1 \leq k \leq M} \tauc(\theta_k).
\end{split}
\end{equation}
In the following we refer to the GCuSum algorithm with $\Theta$ finite as the MCuSum algorithm.
Thus, the GCuSum algorithm \eqref{eq:GLRTCUSUM_algo} has a recursive implementation in this case.\footnote{We note however that while $C_n$ is a non-negative statistic, the statistic $G_n$ can take negative values.} The asymptotic optimality of the GCuSum algorithm, with $\Theta$ finite, is proved in \cite{tart_polu_ISIF_2008}.
Specifically, setting $A=\log M/\alpha$ in \eqref{eq:GLRTCuSum_finitetheta_1} ensures that
\begin{equation}
\label{eq:MCuSumPerf}
\begin{split}
\FAR(\taugc) &\leq \alpha,\\
\WADD^{\theta_k}(\taugc) &\leq \frac{|\log \alpha|}{D(f_{\theta_k} \; ||\; f_0)}(1+o(1)) \mbox{ as } \alpha \to 0, \mbox{ for } 1 \leq k \leq M.
\end{split}
\end{equation}
Thus, due to \eqref{eq:MCuSumPerf}, \eqref{eq:CADDleqWADD} and \eqref{eq:LB},
the GCuSum algorithm is asymptotically optimal for both Problem~\ref{prob:Pollak_Unknownf1} and Problem~\ref{prob:Lorden_Unknownf1}, with $\beta=1$, as $\alpha \to 0$, uniformly over $\theta_k$, $1 \leq k \leq M$.

Now consider the case when the pre- and post-change distributions belong to an exponential family such that
\begin{equation}\label{eq:expo_family_def}
f_\theta(x) = \exp(\theta x-b(\theta)) f_0(x), \theta \in \Theta,
\end{equation}
where, $\Theta$ is an interval on the real line not containing $0$, i.e., $\Theta=[\theta_\ell,\theta_u] \backslash \{0\}$,
and $b(0)=0$.
%Thus, the pre-change distribution corresponds to the parameter taking value $0$.
As claimed in \cite{lord-amstat-1971}, this model can be used to represent a much broader class of
one-parameter exponential family.
%
%by a transformation of parameters, this model
%covers a more general one-parameter exponential family with respect to a sigma finite measure $\nu$,
%where both $f_0$ and $\{f_\theta\}$ belong to that family,
%and the family has been shifted to make the pre-change distribution correspond to parameter value $0$.
%Also, the measure $d\nu$ has been replaced by $e^{b(0)} d\nu$.
For this case, the asymptotic optimality of the GCuSum algorithm is studied in \cite{lord-amstat-1971}.
Specifically, with $\epsilon > 0$, $\Theta'(\alpha) = \{\theta \in [\theta_\ell, \theta_u]: |\theta| > \epsilon\}$ and setting $A=A_\alpha \simeq \log 1/\alpha$ ensures
\begin{equation}
\label{eq:GCuSumPerf_expfam}
\begin{split}
\FAR(\taugc) &\leq \alpha(1+o(1)), \mbox{ as } \alpha \to 0\\
\WADD^{\theta}(\taugc) &\leq \frac{|\log \alpha|}{D(f_{\theta} \; ||\; f_0)}(1+o(1)) \mbox{ as } \alpha \to 0, \mbox{ for all } \theta \in \Theta'(\alpha).
\end{split}
\end{equation}
Here, $\epsilon$ is allowed to decrease to zero as $\alpha \to 0$. As a result, each $\theta \in \Theta$ is covered eventually.
Thus, to detect a change with $\theta$ very close to $0$, we must operate at low false alarm rates.

We remark on the differences between the results in \eqref{eq:MCuSumPerf} and \eqref{eq:GCuSumPerf_expfam}.
While \eqref{eq:MCuSumPerf} is valid only with $\Theta$ finite, the pre- and post-change distributions are allowed to be
arbitrary, and the $\FAR$ result is non-asymptotic. 
On the other hand, in \eqref{eq:GCuSumPerf_expfam}, the distributions are restricted to an exponential
family, and the $\FAR$ result is asymptotic, but the parameter set $\Theta$ is allowed to be uncountably infinite.
We also note that when $\Theta$ is finite, the GCuSum algorithm has a recursive implementation. 

\section{QCD with Observation Control ($\beta<1$), $\theta$ Known}
For the case when $\theta$ is known and $\beta < 1$, in \cite{bane-veer-IT-2013},
we proposed the DECuSum algorithm, which is a two-threshold modification
of the CuSum algorithm \eqref{eq:CUSUM_recursion}, and showed that it is asymptotically optimal
for a variation of both Problem~\ref{prob:Pollak_Unknownf1} and Problem~\ref{prob:Lorden_Unknownf1} (with a different
$\PDC$ metric), for each $\beta$, as $\alpha \to 0$.
Since the duty cycle metric $\PDC$ is different here, in this section we prove the asymptotic
optimality of the DECuSum algorithm with this new definition of the duty cycle metric.

We first describe the DECuSum algorithm. Fix parameters $A \geq 0$, $\mu > 0$, $h \geq 0$. 
Start with $W_0(\theta)=0$. For $n \geq 0$, 
\begin{equation}\label{eq:DECuSumAlgoDesc}
\begin{split}
S_{n+1} &= 1 \text{~only if~} W_{n}(\theta) \geq 0, \\
W_{n+1}(\theta) &= \min\{W_{n}(\theta) + \mu, 0\} \; \text{~if~} S_{n+1} = 0,\\
             &=\left(W_{n}(\theta) + \log \frac{f_{\theta}(X_{n+1})}{f_{0}(X_{n+1})}\right)^{h+} \; \text{~if~} S_{n+1} = 1,\\
\end{split}
\end{equation}
where $(x)^{h+} = \max\{x, -h\}$.

\begin{figure}[htb]
\center
\includegraphics[width=10cm, height=6cm]{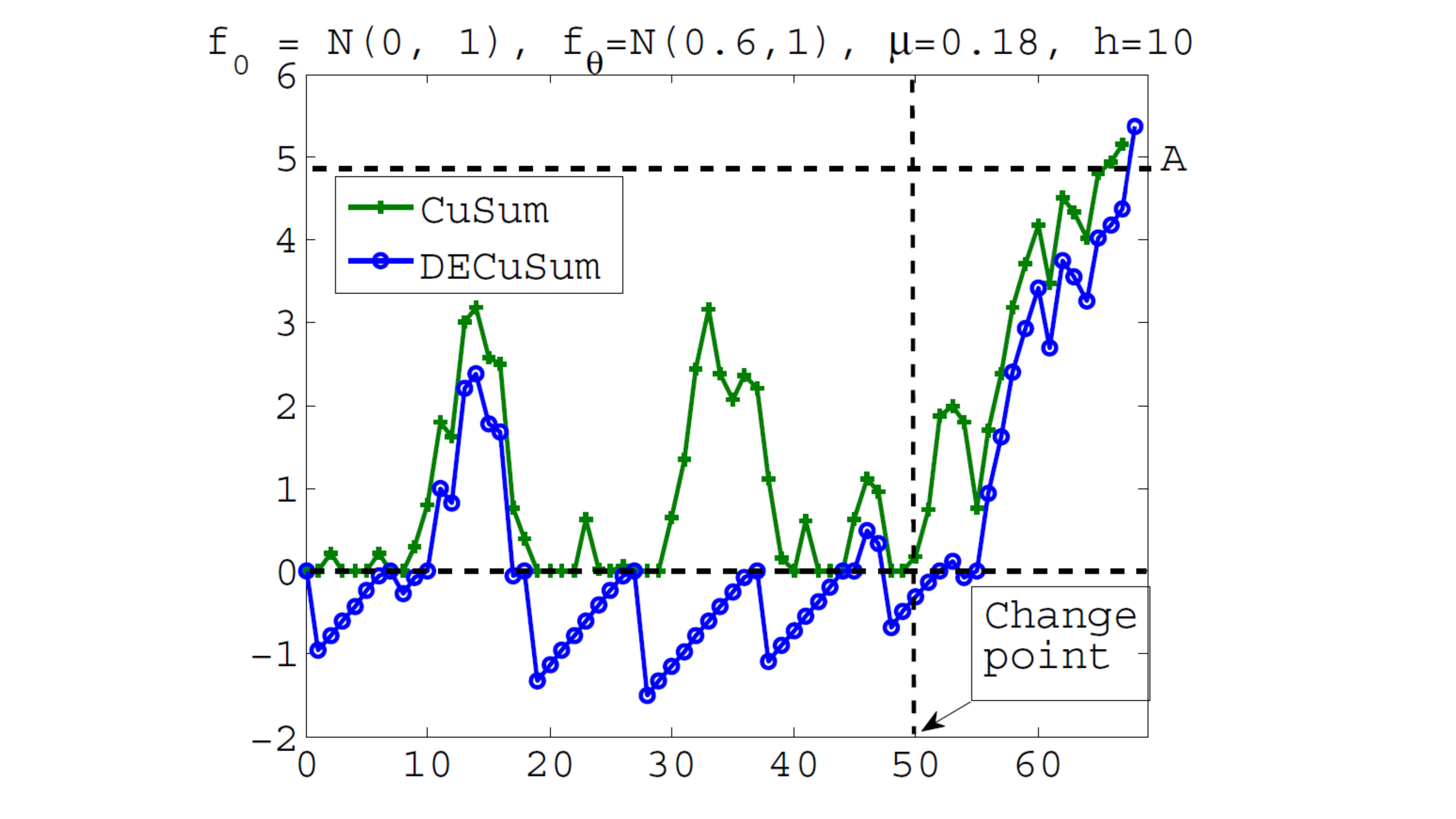}
\caption{Evolution of the CuSum and the DECuSum algorithm for the same set of samples with parameters $f_0=\mathcal{N}(0,1)$, $f_\theta=\mathcal{N}(0.6,1)$, $\mu=0.18$, and $h=10$.}
\label{fig:CuSum_DECuSum_Evo}%\vspace{-0.6cm}
\end{figure}

See Fig.~\ref{fig:CuSum_DECuSum_Evo} for a typical evolution of the CuSum and the DECuSum
algorithms applied to the same set of samples.
When $h=\infty$, the evolution of the DECuSum algorithm can be explained as follows.
Recall that in the CuSum algorithm the log likelihood ratio of the observations is accumulated over time.
If the statistic $C_n(\theta)$ goes below 0, then the statistic is reset to zero. In the DECuSum algorithm,
when the accumulated log likelihood statistic $W_{n}(\theta)$ goes below 0, it is treated as a sign of no change, and
samples are skipped based on the undershoot of the statistic. Mathematically, the statistic is incremented
by a parameter $\mu$ until the statistic reaches 0 from below, at which time the statistic is reset to zero.
This completes a renewal cycle and the above process is repeated till the statistic
$W_{n}(\theta)$ crosses the threshold $A$ from below, at which time a change is declared.
When $h < \infty$, the undershoot of the statistic $W_{n}(\theta)$ is truncated at $-h$, bounding
the maximum number of consecutive samples skipped by $\lceil h/\mu \rceil$. This may be
desired in some applications. The parameters $\mu$ and $h$ are design parameters used
to control the $\PDC$, and the threshold $A$ is used to control the false alarm.

We now prove the asymptotic optimality of the DECuSum algorithm. 
For the theorem below, we need the following definition.
We define the ladder variable \cite{wood-nonlin-ren-th-book-1982}
\[\tau_{-}(\theta) = \inf\left\{n \geq 1: \sum_{k=1}^n \log \frac{f_\theta(X_k)}{f_0(X_k)} < 0\right\}.\]
Then note that $W_{\tau_{-}}(\theta)$ is the ladder height. Recall that $(x)^{h+} = \max\{x, -h\}$.
%In the following theorem, we suppress the dependence on $\theta$ because it is known.
\begin{theorem}
\label{thm:DECuSumOpt}
When the post-change density $f_\theta$ is fixed and known, and $\mu > 0$, $h < \infty$, and $A=|\log \alpha|$, we have
\begin{equation}
\label{eq:DECuSumPerf}
\begin{split}
\FAR(\tauw(\theta)) &\leq \FAR(\tauc(\theta)) \leq \alpha,\\
\PDC(\tauw(\theta)) &=\frac{\Expect_\infty[\tau_{-}(\theta)]}{\Expect_\infty[\tau_{-}(\theta)] + \Expect_\infty[\lceil | W_{\tau_{-}}(\theta)^{h+}|/\mu \rceil]},\\
\WADD^{\theta}(\tauw(\theta)) &\sim \WADD^{\theta}(\tauc(\theta)) \sim \frac{|\log \alpha|}{D(f_\theta \; ||\; f_0)} (1+o(1)) \mbox{ as } \alpha \to 0.
\end{split}
\end{equation}
If $h=\infty$, then
\begin{equation}
\label{eq:PDCApprox}
\PDC(\tauw(\theta)) \leq \frac{\mu}{\mu+D( f_0 \; ||\; f_\theta)}.
\end{equation}
\end{theorem}
\begin{IEEEproof}The proofs for the $\FAR$ and $\WADD$ analysis are identical to that provided in \cite{bane-veer-IT-2013}.
For the $\PDC$ we have the following proof. If $S_n$ is treated as a reward for an
on-off renewal process with the on time distributed according to the law of $\tau_{-}$, and the off time distributed according
to the law of $\lceil |W_{\tau_{-}}|/\mu \rceil$ (with truncation taken into account if $h<\infty$). Then, by the renewal reward theorem we have
\begin{equation*}
\begin{split}
\PDC(\tauw) &= \frac{\Expect_\infty[\tau_{-}]}{\Expect_\infty[\tau_{-}] + \Expect_\infty[\lceil |W_{\tau_{-}}^{h+}/\mu| \rceil]}. %\\
%&\leq \frac{\Expect_\infty[\tau_{-}]}{\Expect_\infty[\tau_{-}] +  \Expect_\infty[W_{\tau_{-}}/\mu]}.
\end{split}
\end{equation*}
This proves \eqref{eq:DECuSumPerf}.

If $h=\infty$, then \eqref{eq:PDCApprox} follows from the above equation
because $x \leq \lceil x \rceil$, and from the Wald's lemma: $\Expect_\infty[|W_{\tau_{-}}|]=\Expect_\infty[\tau_{-}]\; D( f_0 \; ||\; f_\theta)$ \cite{wood-nonlin-ren-th-book-1982}.
\end{IEEEproof}

We note that the expression for the $\PDC$ is not a function of the threshold $A$. Also,
for any given $h>0$, the smaller the value of the parameter $\mu$, the smaller the $\PDC$.

With $A=|\log \alpha|$ and $\mu$ and $h$ set to achieve the $\PDC$ constraint of $\beta$ (independent
of the choice of $A$), the $\WADD$ of the DECuSum algorithm achieves the lower bound \eqref{eq:LB}.
Hence, we have from \eqref{eq:CADDleqWADD} that the algorithm is asymptotically optimal for
both Problem~\ref{prob:Pollak_Unknownf1} and Problem~\ref{prob:Lorden_Unknownf1}, for the given $\beta$, as $\alpha \to 0$.
Thus, the pre-change observation control can be executed, i.e., any arbitrary but fixed fraction of samples
can be dropped before change, without any loss in the asymptotic performance.
%
%We have observed through simulations that the upper bound \eqref{eq:PDCApprox} is a good approximation for the $\PDC$ for large values of $h$.
%We have also shown in \cite{bane-veer-IT-2013} that the DECuSum algorithm performs much
%better that the fractional sampling scheme in which the CuSum algorithm is used and the samples are skipped
%randomly to achieve the constraint on the $\PDC$.

Finally, we note that the DECuSum algorithm can also be described as follows. 
\medskip
\begin{equation}
\label{eq:DECUSUM_AlternateForm}
\begin{split}
\mbox{If } W_{n-1}(\theta)\geq  0, &\\
S_n &=1\\
W_{n}(\theta) &= \max\left\{ -h, \max_{1\leq k \leq n} \; \sum_{i=k}^n \log \frac{f_\theta(X_i^{(S_i)})}{f_0(X_i^{(S_i)})}\right\}.\\
\mbox{If } W_{n-1}(\theta) < 0, &\\
S_n &=0,\\
W_{n}(\theta) &= \min\{0, W_{n-1}(\theta) +  \mu\}.\\
\mbox{Stop at}\quad \quad \quad \;\;&\\
\tauw(\theta) &= \inf\{n \geq 1: W_n(\theta) \geq A\}.
\end{split}
\end{equation}
\medskip
where $\frac{f_\theta(X_i^{(S_i)})}{f_0(X_i^{(S_i)})}=1$ if $S_i=0$.
This description will be useful in Section~\ref{sec:G-DE-CuSum}.

%\vspace{-0.55cm}
\section{The GDECuSum Algorithm}
\label{sec:G-DE-CuSum}
%\vspace{-0.2cm}
In this section we propose the main algorithm of this paper, the GDECuSum algorithm.
This algorithm can be used for the case when
the post-change distribution is composite, and there is a need to perform on-off observation control,
which is the object of study in this paper.
Mathematically, $\beta<1$ in Problem~\ref{prob:Pollak_Unknownf1} and Problem~\ref{prob:Lorden_Unknownf1}, and $\theta$ is unknown.

We now make the important assumption that there exists $\theta^* \in \Theta$
such that $f_{\theta^*}$ is the least favorable distribution among the family $\{f_\theta\}$,
in a sense defined by the following assumption:
%\begin{assumption}
%\[D(f_{\theta^*} \; ||\; f_0) \leq D(f_{\theta} \; ||\; f_0) \quad \forall \theta \in \Theta, \]
%\end{assumption}
%and
\begin{assumption}\label{assum:leastfav}
For each $\theta \in \Theta$,
\[\Expect^{\theta}_1 \left[\log \frac{f_{\theta^*}(X_1)}{f_0(X_1)}\right] = D(f_{\theta} \; ||\; f_0) - D(f_{\theta} \; ||\; f_{\theta^*}) > 0.\]
\end{assumption}

The assumption is satisfied for example when
the law of $\log \frac{f_{\theta^*}(X_1)}{f_0(X_1)}$ under $\{f_{\theta}\}$ is stochastically bounded by
its law under $f_{\theta^*}$ (see Definition 1 in \cite{unni-etal-ieeeit-2011}), i.e.,
\[\Prob^{\theta}_1\left(\log \frac{f_{\theta^*}(X_1)}{f_0(X_1)} > x\right) \geq \Prob^{\theta^*}_1\left(\log \frac{f_{\theta^*}(X_1)}{f_0(X_1)} > x\right), \quad \forall \theta \in \Theta.\]

The latter condition is satisfied for example in the following cases:
\begin{enumerate}
\item $\Theta$ is finite, $\Theta= \{\theta_1, \cdots, \theta_M\}$, $f_0=\mathcal{N}(0,1)$, $f_{\theta_k} = \mathcal{N}(\theta_k,1)$, with
$0<\theta_1 < \theta_2 < \cdots < \theta_M$, and $\theta^*=\theta_1$.
\item $\{f_\theta\}$ and $f_0$ belong to an exponential family such that $f_0=\mathcal{N}(0,1)$, $f_{\theta} = \mathcal{N}(\theta,1)$, with
$\theta \in [0.2, 1]$, and $\theta^*=0.2$.
\end{enumerate}

We now propose the GDECuSum algorithm.
In the GDECuSum algorithm also, just like in the GCuSum algorithm \eqref{eq:GLRTCUSUM_algo},
a family of algorithms are executed in parallel, one for each
post-change parameter, with the difference that the CuSum algorithm corresponding to the parameter
$\theta=\theta^*$ is replaced
by the DECuSum algorithm. Also, the CuSum algorithms corresponding to $\theta \neq \theta^*$ are updated
only when samples are taken.
The least favorable post-change density $f_{\theta^*}$ is used for observation control, 
while the entire family of post-change distributions is used for change detection.

The GDECuSum algorithm is described as follows.
\medskip
\begin{algorithm}
Fix $\mu >0$ and $h \geq 0$,
\begin{equation}
\label{eq:GDECUSUM}
\begin{split}
%V_0(\theta) &= 0  \;\;\; \forall \theta \in \Theta\\
\mbox{Compute for each } & n \geq 1,\\
\bar{G}_{n} &= \max_{1\leq k \leq n} \;\sup_{\theta \in \Theta} \; \sum_{i=k}^n \log \frac{f_\theta(X_i^{(S_i)})}{f_0(X_i^{(S_i)})}. \\
\mbox{If } W_{n-1}(\theta^*) \geq 0,&\\
S_n &=1\\
W_{n}(\theta^*) &= \max\left\{ -h, \max_{1\leq k \leq n} \; \sum_{i=k}^n \log \frac{f_{\theta^*}(X_i^{(S_i)})}{f_0(X_i^{(S_i)})}\right\}.\\
\mbox{If } W_{n-1}(\theta^*) < 0, &\\
S_n &=0 \\
W_{n}(\theta^*) &= \min\{0, W_{n-1}(\theta^*) +  \mu\}.\\
\mbox{Stop at}\quad \quad \quad \quad \;\; &\\
\taugd &= \inf\{n \geq 1: \bar{G}_n \geq A\}.
\end{split}
\end{equation}
\end{algorithm}

The evolution of the GDECuSum algorithm can be described as follows. In this algorithm
two statistics $\bar{G}_n$ and $W_n(\theta^*)$ are computed in parallel. While
the statistic $\bar{G}_n$ is used to detect the change, the statistic $W_n(\theta^*)$
is used for observation control. Specifically, the statistic $W_n(\theta^*)$
is updated using the DECuSum algorithm \eqref{eq:DECUSUM_AlternateForm}.
The statistic $\bar{G}_n$ is
updated using the GCuSum algorithm \eqref{eq:GLRTCUSUM_algo} with the difference that when $W_n(\theta^*)<0$, the
statistic $\bar{G}_n$ is not updated.
%This is because when $W_n(\theta^*)<0$, it is incremented by $\mu$ at each time instant,
%and observations are skipped until $W_n(\theta^*)$ reaches 0 from below. In the absence of
%any new observation, the GCuSum statistic $\bar{G}_n$ cannot be updated.
%In this algorithm, by design, while $W_n(\theta^*)<0$,
%the GCuSum statistic $\bar{G}_n$
%is set to its value in the last time instant.
%This is ensured by the definition $\frac{f_\theta(X_i^{(S_i)})}{f_0(X_i^{(S_i)})}=1$ if $S_i=0$.

Assumption~\ref{assum:leastfav} is critical to the working of this algorithm.
By this assumption
the mean of the log likelihood ratio between $f_{\theta^*}$ and $f_0$ is positive for every possible post-change distribution.
This is because for $\theta \in \Theta$,
\[\Expect^{\theta}_1 \left[\log \frac{f_{\theta^*}(X_1)}{f_0(X_1)}\right] = D(f_{\theta} \; ||\; f_0) - D(f_{\theta} \; ||\; f_{\theta^*}).\]
This ensures that after the change occurs, and after a finite number of samples (irrespective of the threshold $A$),
the DECuSum statistic $W_n(\theta^*)$ always remains positive and no more observations are skipped. This allows
the statistic $\bar{G}_n$ to grow with the right ``slope''.
If the Assumption~\ref{assum:leastfav} is violated, and the post-change parameter is $\theta \neq \theta^*$,
then the statistic $W_n(\theta^*)$ will be below zero for a longer duration of time,
and this time grows to infinity as the threshold $A \to \infty$. Thus, essentially, the growth of the GCuSum statistic
will be intercepted by multiple sojourns of the statistic $W_n(\theta^*)$
below zero. As a result, the change will still be detected, but with a delay larger than the lower bound \eqref{eq:LB}.

\begin{figure}[htb]
\center
\includegraphics[width=10cm, height=6cm]{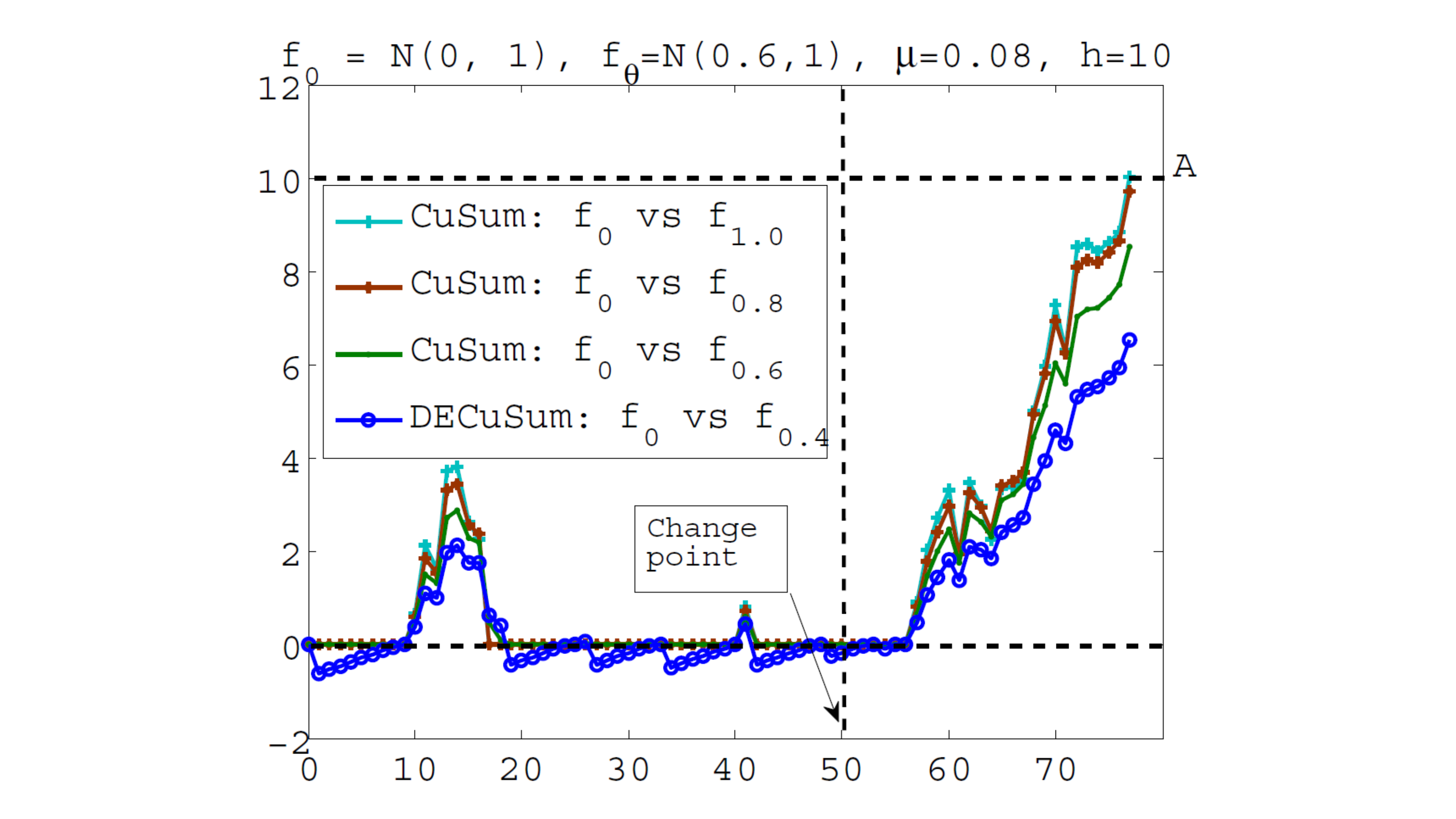}
\caption{Evolution of the GDECuSum algorithm for $f_0=\mathcal{N}(0,1)$, $f_{\theta_1}=\mathcal{N}(0.4,1)$, $f_{\theta_2}=\mathcal{N}(0.6,1)$, $f_{\theta_3}=\mathcal{N}(0.8,1)$, $f_{\theta_4}=\mathcal{N}(1,1)$, $\mu=0.18$, and $h=10$. The post-change parameter is $\theta=\theta_2=0.6$.}
\label{fig:GDECuSum_Evo}%\vspace{-0.6cm}
\end{figure}

For $\Theta=\{\theta_1, \cdots, \theta_M\}$ with $\theta^*=\theta_1$,
the GDECuSum algorithm has a recursive implementation
\footnote{Again note that the statistics $\{\bar{C}_k\}_{k=2}^M$ here are non-negative while $\bar{G}_n$
is allowed to take negative values.}
\begin{equation}
\label{eq:GDECUSUM_Finite_1}
\begin{split}
%V_0(\theta) &= 0  \;\;\; \forall \theta \in \Theta\\
\mbox{If } W_{n-1}(\theta_1) \geq 0, &\\
S_n &= 1,\\
W_{n}(\theta_1) &= \left(W_{n-1}(\theta_1) +  \log \frac{f_{\theta_1}(X_n)}{f_0(X_n)}\right)^{h+}. \\
\mbox{If } W_{n-1}(\theta_1) < 0, &\\
S_n &= 0,\\
W_{n}(\theta_1) &= \min\{0, W_{n-1}(\theta_1) +  \mu\}. \\
\mbox{For } k \geq 2, \quad \quad \quad &\\
\bar{C}_0(\theta_k) &= 0,\\
\bar{C}_{n}(\theta_k) &= \left(\bar{C}_{n-1}(\theta_k) +  \log \frac{f_{\theta_k}(X_n^{(S_n)})}{f_0(X_n^{(S_n)})}\right)^+.\\
\mbox{Stop at}, \quad \quad \quad \quad \; \; &\\
\taugd = \inf\{n \geq 1: & \max\{ W_n(\theta_1), \max_{2\leq k \leq M} \bar{C}_n(\theta_k)\} \geq A\}.
\end{split}
\end{equation}
Thus, for $\Theta$ finite, the GDECuSum algorithm is equivalent to executing $M$ recursive algorithms in parallel.
One is the DE-CuSum algorithm using the least favorable distribution, and the rest $M-1$ algorithms
are the CuSum algorithms. Note that when the DE-CuSum statistic $W_n(\theta_1)<0$, the CuSum statistics
$\{\bar{C}_n(\theta_k)\}_{k=2}^M$ are set to their values in the last time instant.
For the case of finite $\Theta$, we refer to the GDECuSum algorithm by the MDECuSum algorithm.
In Fig.~\ref{fig:GDECuSum_Evo} we plot the evolution of the GDECuSum algorithm (or the MDECuSum algorithm)
for $f_0=\mathcal{N}(0,1)$, $f_{\theta_1}=\mathcal{N}(0.4,1)$, $f_{\theta_2}=\mathcal{N}(0.6,1)$, $f_{\theta_3}=\mathcal{N}(0.8,1)$, $f_{\theta_4}=\mathcal{N}(1,1)$, $\mu=0.18$, and $h=10$. The post-change parameter is $\theta=\theta_2=0.6$.

\section{Asymptotic Optimality of the GDECuSum Algorithm}
\label{sec:GDECuSum_Perf}
The evolution of the GDECuSum algorithm is statistically identical to that of the GCuSum algorithm,
except for the possible sojourns of the statistic $W_n(\theta^*)$ below 0. Also, the sojourn time
of $W_n(\theta^*)$ below zero is completely specified by the DECuSum algorithm.
These two facts will now be used to express the performance of the GDECuSum algorithm
in terms of the performance of the GCuSum algorithm and the DECuSum algorithm.

Define
\[\tauw(\theta^*) = \inf\{n \geq 1: W_n(\theta^*) \geq A\},\]
i.e, $\tauw(\theta^*) $ is the first time the statistic $W_n(\theta^*)$ crosses the threshold $A$.
\begin{theorem}\label{thm:GDECuSum_OPT}
Under the Assumption~\ref{assum:leastfav}, for any fixed $\mu > 0$ and $h\geq 0$ and $A$ we have
\begin{equation}
\label{eq:GDECuSumPerf_FAR_PDC}
\begin{split}
\FAR(\taugd) &\leq \FAR(\taugc),\\
\PDC(\taugd) &= \PDC(\tauw(\theta^*)),
\end{split}
\end{equation}
and for any $\mu > 0$ and $h < \infty$, and any $A \geq 0$,
\begin{equation}
\label{eq:GDECuSumPerf_WADD_secOrd}
\begin{split}
\WADD^{\theta}(\taugd) &\leq \WADD^{\theta}(\taugc) + K_{\mathrm{GD}}, \\
\end{split}
\end{equation}
where $K_{\mathrm{GD}}$ is a constant that is a function of $\mu$ and $h$, but is not a function of $A$.
As a result, for any $\mu > 0$ and $h < \infty$, we have
\begin{equation}
\label{eq:GDECuSumPerf_WADD_FirstOrd}
\begin{split}
\WADD^{\theta}(\taugd) &\sim \WADD^{\theta}(\taugc) \sim \frac{A}{D(f_\theta \; ||\; f_0)}(1+o(1)) \\
& \hspace{1cm} \mbox{ as } A \to \infty, \mbox{ for each } \theta \in \Theta.
\end{split}
\end{equation}
\end{theorem}

We provide the proof of the theorem in the appendix. We now discuss the implications of this result. 
From the theorem we see that, the GDECuSum algorithm can be designed to satisfy any arbitrary $\PDC$ constraint of $\beta$, independent of the choice of $A$.
Also, the $\FAR$ of the GDECuSum algorithm is at least as good as that of the GCuSum algorithm. Finally,
the $\WADD$ of the GDECuSum algorithm is within a constant of the $\WADD$ of the GCuSum algorithm.
From \eqref{eq:CADDleqWADD} and \eqref{eq:CADDWADDWorstGamma_MCuSum} we have
\[
\CADD^\theta(\taugd) \leq \WADD^\theta(\taugd) \leq \WADD^\theta(\taugc) + K_{\mathrm{GD}} = \CADD^\theta(\taugc) + K_{\mathrm{GD}}.
\]

Thus, the GDECuSum algorithm will be asymptotically optimal for the proposed problems for any fixed $\beta$, 
if the GCuSum algorithm is asymptotically optimal for the proposed problems with $\beta=1$. 
%\begin{equation}\label{eq:SetThres_getalpha}
%A=A_\alpha \simeq \log \alpha \mbox{ in \eqref{eq:GLRTCUSUM_algo} } \implies \FAR(\taugc) \leq \alpha,
%\end{equation}
This is formally stated in the next corollary.
\medskip
\begin{corrly}
If the GCuSum algorithm is uniformly asymptotically optimal for a parametric family for Problem~\ref{prob:Pollak_Unknownf1}
or Problem~\ref{prob:Lorden_Unknownf1} with $\beta=1$, 
then under the conditions of the theorem and if $h < \infty$,
the GDECuSum algorithm is also uniformly asymptotically optimal, for the corresponding problem, for each $\beta$, as $\alpha \to 0$.
\end{corrly}
\medskip
Since the GCuSum algorithm is asymptotically optimal (with $\beta=1$) for the two special classes of $\{f_\theta\}$: finite and exponential,
the GDECuSum algorithm is also asymptotically optimal (for each fixed $\beta$) in these two cases. These are stated as corollaries below.

%
%Thus, from \eqref{eq:CADDleqWADD}, setting the threshold $A$ as in the theorem above,
%and choosing the parameters $\mu$ and $h$ to set $\PDC \leq \beta$ ensures
%that the MDECuSum algorithm is asymptotically optimal for both Problem 1 and Problem 2, for the given $\beta$,
%as $\alpha \to 0$. %The following approximation can be used to choose the value of $\mu$:
%%\begin{equation}
%%\label{eq:PDCApprx_MDECuSum}
%%\PDC(\taund) \approx \frac{\mu}{\mu+D( f_0 \; ||\; f_{\theta_1})}.
%%\end{equation}

For a finite family we have the following result.
\medskip
\begin{corrly}
If $\Theta$ if finite, $\Theta = \{\theta_1, \cdots, \theta_M\}$, and Assumption~\ref{assum:leastfav}
is satisfied for some $\theta^* \in \Theta$. Then, for any fixed $\mu > 0$ and $h\geq 0$ and $A=\log M/\alpha$ we have
\begin{equation}
\label{eq:GDECuSumPerf_Fin_FARPDC}
\begin{split}
\FAR(\taugd) &\leq \FAR(\taugc)\leq \alpha,\\
\PDC(\taugd) &= \PDC(\tauw(\theta^*)).
\end{split}
\end{equation}
Also, if $\mu > 0$ and $h < \infty$, then
\begin{equation}
\label{eq:GDECuSumPerf_Fin_WADD}
\begin{split}
\WADD^{\theta}(\taugd) &\sim \WADD^{\theta}(\taugc) \sim \frac{|\log \alpha|}{D(f_{\theta_k} \; ||\; f_0)}(1+o(1)) \\
&  \mbox{ as } \alpha \to 0, \mbox{ for each } \theta_k, k=1,\cdots,M.
\end{split}
\end{equation}
\end{corrly}
\begin{IEEEproof}
The result follows from \eqref{eq:MCuSumPerf} and Theorem~\ref{thm:GDECuSum_OPT}.
\end{IEEEproof}
\medskip
For one-parameter exponential families we have the following result.
\medskip
\begin{corrly}
If $\{f_\theta\}$, $f_0$ belong to a one-parameter exponential family, i.e., if the following is satisfied,
\[
f_\theta(x) = \exp(\theta x-b(\theta)) f_0(x), \mbox{ for } \theta \in \Theta,
\]
where, $\Theta=[\theta_\ell,\theta_u]$, with $0 < \theta_\ell < \theta_u$, and $b(0)=0$.
Also, Assumption~\ref{assum:leastfav} is satisfied for some $\theta^*\in\Theta$.
Then, for any fixed $\mu > 0$, $h\geq 0$ and $A=A_\alpha \simeq \log 1/\alpha$ we have
\begin{equation}
\label{eq:GDECuSumPerf_Expo_FARPDC}
\begin{split}
\FAR(\taugd) &\leq \FAR(\taugc)\leq \alpha(1+o(1)), \mbox{ as } \alpha \to 0\\
\PDC(\taugd) &= \PDC(\tauw(\theta^*)).
\end{split}
\end{equation}
And if $h < \infty$, then
\begin{equation}
\label{eq:GDECuSumPerf_Expo_WADD}
\begin{split}
\WADD^{\theta}(\taugd) &\sim \WADD^{\theta}(\taugc) \sim \frac{|\log \alpha|}{D(f_\theta \; ||\; f_0)}(1+o(1)) \\
&  \mbox{ as } \alpha \to 0, \mbox{ for each } \theta \in [\theta_\ell, \theta_u].
\end{split}
\end{equation}
\end{corrly}

\begin{IEEEproof}
The result follows from \eqref{eq:GCuSumPerf_expfam}  and Theorem~\ref{thm:GDECuSum_OPT}.
\end{IEEEproof}
\medskip

Since, the GDECuSum algorithm achieves the lower bound \eqref{eq:LB}, the algorithm is asymptotically optimal for the two cases
specified in the corollaries above, for a given $\beta$, uniformly over $\theta \in \Theta$, as $\alpha \to 0$.
%\medskip

\section{Discussion}
\label{sec:Discussion}
In this section we discuss possible extensions of the results developed in the previous sections.

\subsection{Extension to Mixture Based Tests}
\label{sec:Mixture}
In the classical QCD problem with unknown post-change distribution, an alternative to 
the GLRT based approach is a mixture based approach. Specifically, let $\pi(\theta)$ be a probability measure on the 
parameter space $\Theta$. Then, a mixture based CuSum test is given by 
\begin{equation}
\label{eq:MixCUSUM_algo}
\begin{split}
\tilde{G}_{n} &= \max_{1\leq k \leq n} \; \log \int_{\theta \in \Theta} \; \prod_{i=k}^n \frac{f_\theta(X_i)}{f_0(X_i)} d\pi(\theta). \quad \mbox{ for } n\geq 1,\\
\taugc &= \inf\{n \geq 1: \tilde{G}_n \geq A\},
\end{split}
\end{equation}
It is well known that under some conditions the mixture based test is also uniformly asymptotically optimal, 
for both Problem~\ref{prob:Pollak_Unknownf1} and Problem~\ref{prob:Lorden_Unknownf1} with $\beta=1$, as $\alpha \to 0$; 
see \cite{tart-niki-bass-2014}. 

Similar to the GDECuSum algorithm, one can also define a data-efficient extension of the above mixture based CuSum test, 
when a least favourable member is present in the post-change family of distributions. 
The preceeding analysis on the GDECuSum algorithm will also hold true almost verbatim for the mixture based data-efficient test, 
with the expection of the argument of type provided in \eqref{eq:ProbSuccLB}. 
For the mixture based data-efficient test \eqref{eq:ProbSuccLB} has to be replaced by the following equation:  
\begin{equation}
\begin{split}
\Prob^\theta_1 &(\mbox{success in $1^{st}$ cycle}) \\
&=\Prob^\theta_1(\mbox{Statistic $\tilde{G}_n $ reaches $A$ before $\sum_{k=1}^n \log \frac{f_{\theta^*}(X_k)}{f_0(X_k)}$ goes below $-w$})\\
&\geq \Prob^\theta_1\left(\sum_{k=1}^n \log \frac{f_{\theta^*}(X_k)}{f_0(X_k)} \geq 0, \; \forall n\right)=q_\theta.
\end{split}
\end{equation}
The above equation will be valid if the mixture distribution $\pi(\theta)$ is chosen such that for each $\theta \in \Theta$, 
$\tilde{G}_n \to \infty$ a.s. $\Prob^\theta_1$. 

\subsection{Extension to Window Limited Tests}
Recall that unless the post-change distribution belongs to a finite family, the GDECuSum algorithm does not have a recursive implementation.
In the classical QCD literature, this problem is addressed by proposing window based tests; see Lai \cite{lai-ieeetit-1998}. 
It is straighforward to show that the data-efficient extensions of such window based GLRT and mixture based tests 
also retain the asymptotic optimality properties of the GDECuSum algorithm.

\subsection{Extension to Parametric Families with no Least Favorable Distribution}
\label{sec:ExistOfLeastFavDist}
One fundamental assumption in this paper is the existence of a least favourable distribution in the post-change family 
 in the sense of Assumption~\ref{assum:leastfav}, i.e., there is a distribution
$f_{\theta^*}$ such that for each $\theta \in \Theta$,
\[
\Expect^{\theta}_1 \left[\log \frac{f_{\theta^*}(X_1)}{f_0(X_1)}\right] > 0.
\]
We used this assumption in the proof of Theorem~\ref{thm:GDECuSum_OPT}.  
The positive mean of the log likelihood ratio $\log \frac{f_{\theta^*}(X_1)}{f_0(X_1)}$ under each $\theta$ ensures
that after a finite number of time slots, no observations are skipped using the DE-CuSum algorithm, and the change is detected efficiently.

However, for a given parametric family, there may not be a distribution that satisfies Assumption~\ref{assum:leastfav}. In such a case,
the results of this paper can be extended to cases where a distribution $g$ exists satisfying the assumption,
i.e.,
\[
\Expect^{\theta}_1 \left[\log \frac{g(X_1)}{f_0(X_1)}\right] > 0, \; \forall \theta \in \Theta.
\]
Thus, as long as such a distribution exists, we can design the DE-CuSum algorithm using the distribution $g$ and the positive drift
in the last equation will ensure that the GDECuSum with this new modification is still asymptotically optimal.
We however note that in the proof of Theorem~\ref{thm:GDECuSum_OPT} we used the fact that $\theta^* \in \Theta$. Since $g$ may not be in the
parametric family, the proof needs to be modified. This can be accomplished by replacing the arguments in \eqref{eq:ProbSuccLB} with
\begin{equation}
\begin{split}
\Prob^\theta_1 &(\mbox{success in $1^{st}$ cycle}) = \Prob^\theta_1(G_{\tau_1(w)}>A)\\
&=\Prob^\theta_1(\mbox{Statistic $G_n $ reaches $A$ before $\sum_{k=1}^n \log \frac{g(X_k)}{f_0(X_k)}$ goes below $-w$})\\
&=\Prob^\theta_1\left(\mbox{$\max_{1\leq k \leq n} \;\sup_{\theta \in \Theta} \; \sum_{i=k}^n \log \frac{f_\theta(X_i)}{f_0(X_i)}$ reaches $A$}\right.\\
&\hspace{2cm}\left.\mbox{before $\sum_{k=1}^n \log \frac{g(X_k)}{f_0(X_k)}$ goes below $-w$}\right)\\
&\geq \Prob^\theta_1\left(\sum_{k=1}^n \log \frac{g(X_k)}{f_0(X_k)} \geq 0, \; \forall n\right).
\end{split}
\end{equation}
The last quantity is positive because $\Expect^{\theta}_1 \left[\log \frac{g(X_1)}{f_0(X_1)}\right] > 0$ \cite{wood-nonlin-ren-th-book-1982}.

\section{Numerical Results}
\label{sec:DECompositeNumerical}

\begin{figure}[htb]
\center
%\vspace{-0.3cm}
\includegraphics[width=10cm, height=6cm]{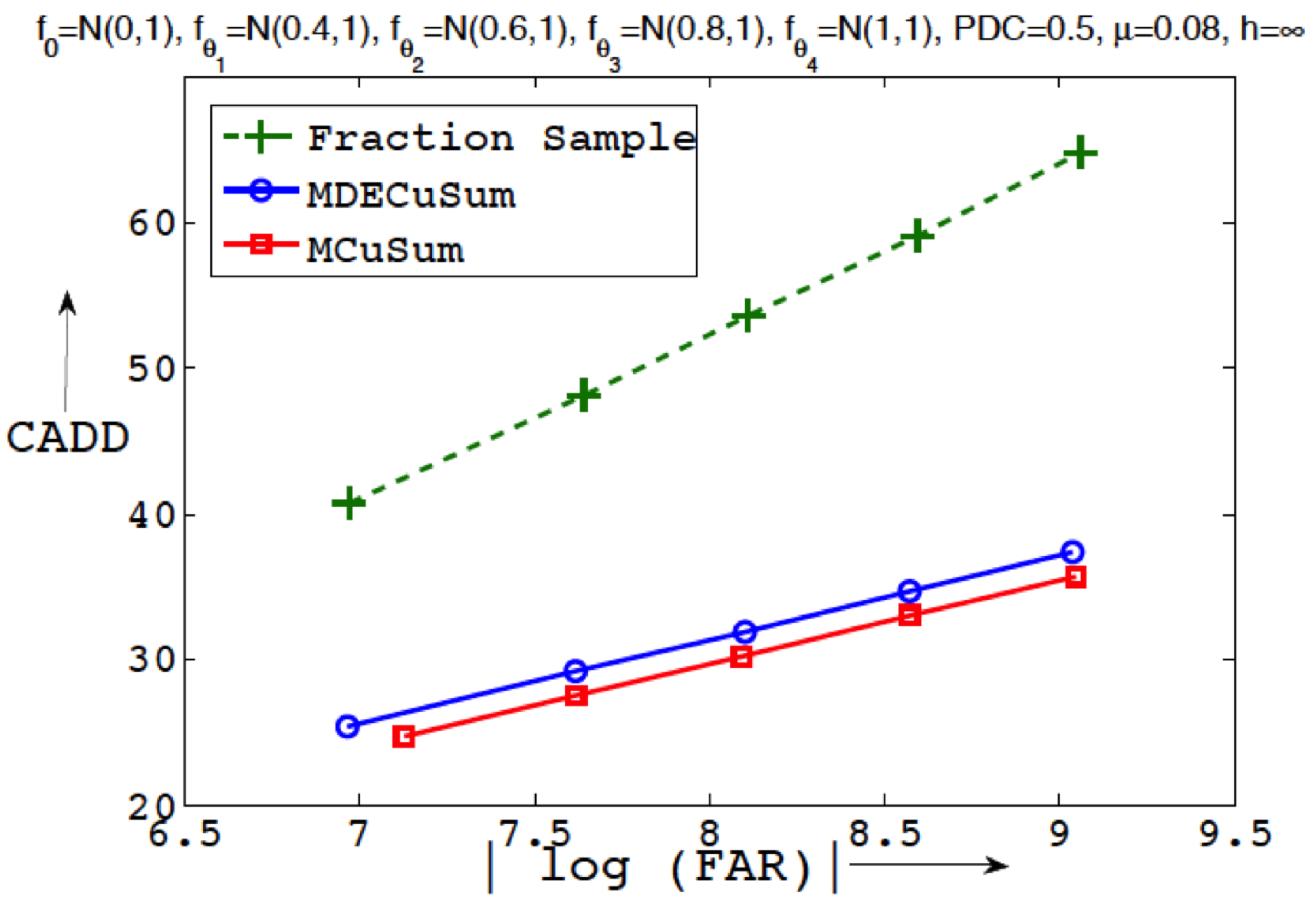}
\caption{Comparative performance of the GDECuSum algorithm, the GCuSum algorithm,
and the fractional sampling scheme. The post-change parameter is $\theta=\theta_2=0.6$.}
\label{fig:NDE_NC_Frac}%\vspace{-0.6cm}
\end{figure}

In Fig.~\ref{fig:NDE_NC_Frac} we plot the $\CADD$--$\FAR$ trade-off curves obtained using simulations
for the
GDECuSum algorithm \eqref{eq:GDECUSUM}, the GCuSum algorithm \eqref{eq:GLRTCUSUM_algo},
and the fractional sampling scheme. In the latter, the GCuSum algorithm is used
and observations are skipped
randomly, independent of the observation process.
The simulation set used is: $M=4$, $f_0=\mathcal{N}(0,1)$, $f_{\theta_1}=\mathcal{N}(0.4,1)$,
$f_{\theta_2}=\mathcal{N}(0.6,1)$, $f_{\theta_3}=\mathcal{N}(0.8,1)$, $f_{\theta_4}=\mathcal{N}(1,1)$,
$\mu=0.08$ and $h=\infty$. The post-change parameter is $\theta=\theta_2=0.6$,
and the value of $\mu$ is chosen using \eqref{eq:PDCApprox}
and \eqref{eq:GDECuSumPerf_FAR_PDC}
 to achieve a $\PDC=0.5$ (skipping/saving 50\% of the samples).
To achieve a $\PDC$ of $0.5$ through the fractional sampling scheme, every alternate sample
is skipped in the GCuSum algorithm.
In the figure we see that skipping samples randomly results
in a twofold increase in delay as compared to that of the GCuSum algorithm.
However, if we use the GDECuSum algorithm and use the state of the system to skip observations,
then there is a small and constant penalty on the delay, as compared to the performance
of the GCuSum algorithm. Thus, the GDECuSum algorithm provides a significant gain in performance as
compared to the fractional sampling scheme.

\section{Conclusions}
\label{sec:Conclusions}
\vspace{-0.1cm}
We have extended our work on data-efficient quickest change detection in \cite{bane-veer-IT-2013}
to the case when the post-change distribution is composite. 
%We have shown that the right recipe is a combination of GLRT and robust approaches. 
If the post-change family of distribution has a least favourable member, 
then we have proposed an algorithm in which, the observation control is implemented 
using the least favorable member, and the change is detected using a GLRT based approach. 
We have shown that under standard conditions used in the literature, the 
proposed algorithm is asymptotically optimal. The implication is that 
an arbitrary but fixed fraction 
of observations can be skipped before change, without affecting the asymptotic performance, 
and this can be done even when the post-change distribution is composite. 
Our numerical results for moderate values of false alarm rate
show that the GDECuSum algorithm incurs a small bounded delay penalty 
relative to the GCuSum algorithm, for considerable gains in data-efficiency. 
This is in sharp contrast to the gain obtained by the commonly used approach to data/energy efficiency based on fractional sampling.  
We have also shown that this work can be extended to prove optimality of data-efficient extensions of mixture based tests 
and window limited tests. Furthermore we shown that as long as there is a distribution that 
is not necessarily part of the family, but can serve the role of a least favourable distribution,  
a data-efficient test can be designed that is asymptotically optimal.
%An interesting direction for future research is to introduce the concept of data-efficiency 
%in more general QCD models. 

\section*{Appendix}
\begin{IEEEproof}[Proof of Theorem~\ref{thm:GDECuSum_OPT}]
We recall that the GCuSum algorithm is the GLRT based test discussed in \eqref{eq:GLRTCUSUM_algo},
and the GDECuSum algorithm is its data-efficient modification discussed in \eqref{eq:GDECUSUM}, where the
observation control is executed based on the least favorable distribution $f_{\theta^*}$.

We wish to prove \eqref{eq:GDECuSumPerf_FAR_PDC} and \eqref{eq:GDECuSumPerf_WADD_secOrd}, i.e.,
for any $\mu>0$, $h\geq 0$, and $A$,
\begin{equation*}
\begin{split}
\FAR(\taugd) &\leq \FAR(\taugc),\\
\PDC(\taugd) &= \PDC(\tauw(\theta^*)),
\end{split}
\end{equation*}
and for any $\mu > 0$ and $h < \infty$, and any $A$,
\begin{equation*}
\begin{split}
\WADD^{\theta}(\taugd) &\leq \WADD^{\theta}(\taugc) + K_{\text{GD}}, \\
\end{split}
\end{equation*}
where $K_{\text{GD}}$ is a constant that is a function of $\mu$ and $h$, but is not a function of $A$.

The $\PDC$ result follows from the $\PDC$ result proved in
Theorem~\ref{thm:DECuSumOpt} because the observation control is governed by the statistic $W_n(\theta^*)$.
We now prove the $\FAR$ and the $\WADD$ results. Both the results are based on the idea that
the evolution of the GDECuSum algorithm is statistically identical to that of the GCuSum algorithm $\taugc$,
except for the possible sojourns of the statistic $W_n(\theta^*)$ below 0. 

Under $\Prob_\infty$, because of the i.i.d. nature of the observations, the sojourns
of the statistic $W_n(\theta^*)$ below $0$ only leads to a larger mean time to false alarm for the GDECuSum algorithm.

On the other hand, under each $\Prob^\theta_1$, the average number of times the statistic $W_n(\theta^*)$ goes below $0$
is bounded by a constant, which is not a function of $A$. This is due to the fact that $f_{\theta^*}$ is the least
favorable distribution, and as a result the drift of $W_n(\theta^*)$ is positive. Since $h < \infty$, the mean time spent by the statistic $W_n(\theta^*)$
each time it goes below $0$, it bounded by $\lceil h/\mu \rceil$. Thus, the total average mean time spent
by the statistic $W_n(\theta^*)$ below $0$ is bounded above by a constant.
This in turn guarantees that the delay of the GDECuSum algorithm is within a constant
of the GCuSum algorithm. The rest of the proof below formalizes these arguments.

We start by writing the stopping time $\taugc$ as a sum of a random number of stopping times. Such a
representation is critical to this proof.
Toward this end we define a set of
new stopping variables. 
Let $w \in [0, A)$, and define
\[
\tau_1(w) = \inf\left\{n\geq 1: G_n > A \mbox{ or } \; w + \sum_{k=1}^n \log \frac{f_{\theta^*}(X_k)}{f_0(X_k)} < 0\right\}.
\]
This is the first time for either the GCuSum statistic $G_n$ to hit $A$ or the random walk $\sum_{k=1}^n \log \frac{f_{\theta^*}(X_k)}{f_0(X_k)}$ to go below $-w$.

On paths over which $G_{\tau_1(w)}<A$, let
\[
\tau_2(w) = \inf\left\{n> \tau_1(w): G_n > A \mbox{ or } \sum_{k=\tau_1(w)+1}^n \log \frac{f_{\theta^*}(X_k)}{f_0(X_k)} < 0\right\}.
\]
Thus, on paths such that $G_{\tau_1(w)}<A$,
after the time $\tau_1(w)$,
the time $\tau_2(w)$ is the first time for either the statistic $G_n$ to cross $A$ or the
random walk $\sum_{k=\tau_1(w)+1}^n \log \frac{f_{\theta^*}(X_k)}{f_0(X_k)}$ to go below $0$.
We define, $\tau_3(w)$, etc. similarly. Next let,
\[
N(w) = \inf\{k \geq 1: G_{\tau_k} > A\}.
\]
For simplicity we introduce the notion of ``cycles'', ``success'' and ``failure''. With reference to the
definitions of $\tau_k(w)$'s above, we say that a success has occurred if the statistic $G_n$
crosses $A$ before the random walk $\sum_{k=1}^n \log \frac{f_{\theta^*}(X_k)}{f_0(X_k)}$ goes below $-w$.
In that case we also say that the number of cycles to $A$ is 1. If on the other hand,
the random walk $\sum_{k=1}^n \log \frac{f_{\theta^*}(X_k)}{f_0(X_k)}$ goes below $-w$ before $G_n$ crosses $A$,
we say a failure has occurred.
The number of cycles is 2, if now
the statistic $G_n$ crosses $A$ before the random walk
$\sum_{k=\tau_1(w)+1}^n \log \frac{f_{\theta^*}(X_k)}{f_0(X_k)}$ goes below $0$.
Thus,
$N(w)$ is the number of cycles to success.

We note that for any given $\theta$,
\[
N(w) \leq \taugc \leq \tauc(\theta).
\]
This is because each cycle has length at least $1$, and $\tauc(\theta)$ is nothing but the $\taugc$ without the
$\sup$ over $\Theta$.
Since, $\tauc(\theta)$ is finite a.s. under both $\Prob_\infty$ and $\Prob^\theta_1$, for each $\theta \in \Theta$ (see Lorden \cite{lord-amstat-1971}), even
$N(w) < \infty$ a.s. under both $\Prob_\infty$ and $\Prob^\theta_1$, for any $\theta \in \Theta$.

Define $\lambda_1(w) = \tau_1(w)$, $\lambda_2(w)= \tau_2(w) - \tau_1(w)$, etc.
Then we in fact have
\begin{equation}\label{eq:GCuSum_sumoflambda}
\taugc = \sum_{k=1}^{N(w)}\lambda_k(w).
\end{equation}
An important point to observe here is that while the terms on the right-hand side
depend on $w$, their sum does not and equals $\taugc$.

We now bound the mean of $N(w)$ under $\Prob^\theta_1$ by a number that is not a function of $w$ and
threshold $A$.
With the identity
\[
\Expect^\theta_1[N(w)] = \sum_{k=1}^\infty \Prob^\theta_1(N(w) \geq k)
\]
in mind, and using the terminology of cycles, success and failure just defined, we write
\begin{equation*}
\begin{split}
\Prob^\theta_1(N(w) \geq k) &= \Prob^\theta_1(\mbox{fail in 1st cycle}) \\
                   &\; \; \cdots \Prob^\theta_1(\mbox{fail in $k-1^{st}$ cycle}|\mbox{fail in all previous}).
\end{split}
\end{equation*}
Now,
\begin{equation*}
\begin{split}
\Prob^\theta_1(\mbox{fail in $i^{th}$ cycle}|&\mbox{fail in all previous}) \\
&= 1 - \Prob^\theta_1(\mbox{success in $i^{th}$ cycle}|\mbox{fail in all previous}).
\end{split}
\end{equation*}

We claim that
\begin{equation}\label{eq:Claim_ThmGDECuSumOpt}
\Prob^\theta_1(\mbox{success in $i^{th}$ cycle}|\mbox{fail in all previous}) \geq \Prob^\theta_1\left(\sum_{k=1}^n \log \frac{f_{\theta^*}(X_k)}{f_0(X_k)} \geq 0, \; \forall n \right).
\end{equation}
From \cite{wood-nonlin-ren-th-book-1982} it is well known that $\Prob^\theta_1(\sum_{k=1}^n \log \frac{f_{\theta^*}(X_k)}{f_0(X_k)} \geq 0, \; \forall n) > 0$. This is because under $\Prob^\theta_1$, by the Assumption~\ref{assum:leastfav},
the drift of the random walk $\sum_{k=1}^n \log \frac{f_{\theta^*}(X_k)}{f_0(X_k)}$ is positive.
Thus, if
\[
q_\theta = \Prob^\theta_1 \left(\sum_{k=1}^n \log \frac{f_{\theta^*}(X_k)}{f_0(X_k)} \geq 0, \; \forall n \right),
\]
then,
\[
\Prob^\theta_1(N(w) \geq k) \leq (1-q_\theta)^{k-1}.
\]
Note that the right-hand side is not a function of the initial point $w$, nor is a function of the threshold $A$.
Hence,
\begin{equation}\label{eq:ENx_UB_11}
\Expect^\theta_1[N(w)] = \sum_{k=1}^\infty \Prob^\theta_1(N(w) \geq k) \leq \sum_{k=1}^\infty (1-q_\theta)^{k-1} = \frac{1}{q_\theta} < \infty.
\end{equation}

To prove the above claim \eqref{eq:Claim_ThmGDECuSumOpt} we note that
\begin{equation}\label{eq:ProbSuccLB}
\begin{split}
\Prob^\theta_1 &(\mbox{success in $1^{st}$ cycle}) = \Prob^\theta_1(G_{\tau_1(w)}>A)\\
&=\Prob^\theta_1(\mbox{Statistic $G_n $ reaches $A$ before $\sum_{k=1}^n \log \frac{f_{\theta^*}(X_k)}{f_0(X_k)}$ goes below $-w$})\\
&=\Prob^\theta_1\left(\mbox{$\max_{1\leq k \leq n} \;\sup_{\theta \in \Theta} \; \sum_{i=k}^n \log \frac{f_\theta(X_i)}{f_0(X_i)}$ reaches $A$}\right.\\
&\hspace{2cm}\left.\mbox{before $\sum_{k=1}^n \log \frac{f_{\theta^*}(X_k)}{f_0(X_k)}$ goes below $-w$}\right)\\
&\geq\Prob^\theta_1\left(\mbox{$\max_{1\leq k \leq n} \;\sum_{i=k}^n \log \frac{f_{\theta^*}(X_i)}{f_0(X_i)}$ reaches $A$}\right.\\
&\hspace{2cm}\left.\mbox{before $\sum_{k=1}^n \log \frac{f_{\theta^*}(X_k)}{f_0(X_k)}$ goes below $-w$}\right) \\
&\geq\Prob^\theta_1\left(\mbox{$\sum_{k=1}^n \log \frac{f_{\theta^*}(X_k)}{f_0(X_k)}$ reaches $A$}\right.\\
&\hspace{2cm}\left.\mbox{before $\sum_{k=1}^n \log \frac{f_{\theta^*}(X_k)}{f_0(X_k)}$ goes below $-w$}\right) \\
&\geq \Prob^\theta_1\left(\sum_{k=1}^n \log \frac{f_{\theta^*}(X_k)}{f_0(X_k)} \geq 0, \; \forall n\right)=q_\theta.
\end{split}
\end{equation}
Here, the first inequality follows because $\theta^* \in \Theta$ over which the supremum is being taken. The last inequality follows
because $\sum_{k=1}^n \log \frac{f_{\theta^*}(X_k)}{f_0(X_k)} \to \infty$ a.s. under $\Prob^\theta_1$ since $\theta^*$ is least favorable.

For the second cycle note that
\begin{equation*}
\begin{split}
\Prob^\theta_1 &(\mbox{success in $2^{nd}$ cycle} | \mbox{failure in first})
 = \Prob^\theta_1 \left(G_{\tau_2(w)}>A|G_{\tau_1(w)}<A) \right)\\
&=\Prob^\theta_1\left(\mbox{Statistic $G_n$, $n> \tau_1(w)$, reaches $A$} \right.\\
   &\hspace{2cm} \left. \mbox{ before $\sum_{k=\tau_1(w)+1}^n \log \frac{f_{\theta^*}(X_k)}{f_0(X_k)}$ goes below $0$} \; \big| \;
     G_{\tau_1(w)}<A\right)\\
&=\Prob^\theta_1\left(\mbox{$\max_{1\leq k \leq n} \;\sup_{\theta \in \Theta} \; \sum_{i=k}^n \log \frac{f_\theta(X_i)}{f_0(X_i)}$,
    $n> \tau_1(w)$, reaches $A$}\right.\\
   &\hspace{2cm}\left. \mbox{before $\sum_{k=\tau_1(w)+1}^n \log \frac{f_{\theta^*}(X_k)}{f_0(X_k)}$ goes below $0$} \;\big| \;
    G_{\tau_1(w)}<A\right)\\
&\geq\Prob^\theta_1\left(\mbox{$\max_{\tau_1(w) < k \leq n} \;\sum_{i=k}^n \log \frac{f_{\theta^*}(X_i)}{f_0(X_i)}$, for $n>
    \tau_1(w)$, reaches $A$}\right.\\
   &\hspace{2cm}\left. \mbox{before $\sum_{k=\tau_1(w)+1}^n \log \frac{f_{\theta^*}(X_k)}{f_0(X_k)}$ goes below $0$}|
    G_{\tau_1(w)}<A\right) \\
&\geq\Prob^\theta_1\left(\mbox{$\sum_{k=\tau_1(w)+1}^n \log \frac{f_{\theta^*}(X_k)}{f_0(X_k)}$, for $n> \tau_1(w)$, reaches
    $A$}\right.\\
   &\hspace{2cm}\left.\mbox{before $\sum_{k=\tau_1(w)+1}^n \log \frac{f_{\theta^*}(X_k)}{f_0(X_k)}$ goes below $0$}|
    G_{\tau_1(w)}<A\right) \\
&=\Prob^\theta_1\left(\mbox{$\sum_{k=1}^n \log \frac{f_{\theta^*}(X_k)}{f_0(X_k)}$ reaches $A$}\right.\\
   &\hspace{2cm}\left.\mbox{before $\sum_{k=1}^n \log \frac{f_{\theta^*}(X_k)}{f_0(X_k)}$ goes below $0$}\right) \\
&\geq \Prob^\theta_1\left(\sum_{k=1}^n \log \frac{f_{\theta^*}(X_k)}{f_0(X_k)} \geq 0, \; \forall n \right)=q_\theta.
\end{split}
\end{equation*}
Almost identical arguments for the other cycles proves the claim that
\[
\Prob^\theta_1(\mbox{success in $i^{th}$ cycle}|\mbox{fail in all previous}) \geq \Prob^\theta_1\left(\sum_{k=1}^n \log \frac{f_{\theta^*}(X_k)}{f_0(X_k)} \geq 0, \; \forall n\right),
\]
and hence it follows that
\[
\Expect^\theta_1[N(w)] \leq\frac{1}{q_\theta} < \infty.
\]

Let
\[
\taugd(w) = \inf\{n \geq 1: \bar{G}_n > A, \mbox{ with } W_0(\theta^*)=w\}.
\]
Clearly, $\taugd=\taugd(0)$.

Just like we did for $\taugc$, we now write the time $\taugd(w)$ as a sum of stopping times.
We will then draw parallels between representation
of this type for $\taugc$ and $\taugd(w)$ to prove the theorem.

Note that the sojourn of the statistic $\bar{G}_n$ to $A$ may include
alternate sojourns of the statistic $W_n(\theta^*)$ above and below $0$.
Motivated by this we define a set of
new variables.
Let $w \in [0, A)$, and define
\[
\bar{\tau}_1(w) = \inf\left\{n\geq 1: \bar{G}_n > A \mbox{ or } \; W_n(\theta^*) < 0; \mbox{ starting with } W_0(\theta^*)=w \right\}.
\]
This is the first time for either the GDECuSum statistic $\bar{G}_n$ to hit $A$ or the DE-CuSum statistic
$W_n(\theta^*)$ to go below $0$, starting with $W_0(\theta^*)=w$.
On paths over which $\bar{G}_{\bar{\tau}_1(w)}<A$, let
$t_1(w)$ be the number of consecutive samples skipped after $\bar{\tau}_1(w)$ using the DE-CuSum statistic. On
such paths again, let
\[
\bar{\tau}_2(w) = \inf\left\{n> \bar{\tau}_1(w)+t_1(w): \bar{G}_n > A \mbox{ or } W_n(\theta^*) < 0 \right\}.
\]
Thus, on paths such that $\bar{G}_{\bar{\tau}_1(w)}<A$,
after the time $\bar{\tau}_1(w)+t_1(w)$,
the time $\bar{\tau}_2(w)$ is the first time for $\bar{G}_n$ to either cross $A$ or the DE-CuSum statistic
$W_n(\theta^*)$ to go below $0$.
We define, $t_2(w)$, $\bar{\tau}_3(w)$, etc. similarly. Next let
\[
\bar{N}(w) = \inf\{n \geq 1: \bar{G}_{\bar{\tau}_n} > A\}.
\]
We also define
$\bar{\lambda}_1(w) = \bar{\tau}_1(w)$, $\bar{\lambda}_2(w)= \bar{\tau}_2(w) - \bar{\tau}_1(w) -t_1(w)$, etc.

We now make an important observation. We observe that due to the i.i.d. nature of the observations
\begin{equation}
\label{eq:N_Lambda_same}
\begin{split}
\bar{N}(w) &\stackrel{d}{=} N(w)\\
\bar{\lambda}_k(w) &\stackrel{d}{=} \lambda_k(w), \; \forall k.
\end{split}
\end{equation}
In fact, we also have, 
\begin{equation}\label{eq:sum_lambdas_equal}
\sum_{n=1}^{\bar{N}(w)}\bar{\lambda}_k(w) \stackrel{d}{=} \sum_{n=1}^{N(w)} \lambda_k(w), \; \forall w.
\end{equation}
Then we have
\[
\bar{N}(w) < \infty \mbox{ a.s. under both $\Prob_\infty$ and $\Prob^\theta_1$ for each $\theta\in\Theta$},
\]
and
\[
\taugd(w) = \sum_{k=1}^{\bar{N}(w)}\bar{\lambda}_k(w) + \sum_{k=1}^{\bar{N}(w)-1}t_k(w).
\]

We are now ready to prove the $\FAR$ result. Using \eqref{eq:N_Lambda_same}, \eqref{eq:sum_lambdas_equal} and \eqref{eq:GCuSum_sumoflambda},
and the observation following \eqref{eq:GCuSum_sumoflambda}, we have
\begin{equation}\label{eq:GDECuSum_FARProofStmnt}
\begin{split}
\Expect_\infty[\taugd] = \Expect_\infty[\taugd(0)]
&=\Expect_\infty\left[\sum_{n=1}^{\bar{N}(0)}\bar{\lambda}_k(0)\right] + \Expect_\infty\left[\sum_{n=1}^{\bar{N}(0)-1}t_k(0)\right]\\
&=\Expect_\infty\left[\sum_{n=1}^{N(0)} \lambda_k(0)\right] + \Expect_\infty\left[\sum_{n=1}^{N(0)-1} t_k(0)\right]\\
&=\Expect_\infty\left[\taugc \right] + \Expect_\infty\left[\sum_{n=1}^{N(0)-1} t_k(0)\right]\\
&\geq \Expect_\infty\left[\taugc \right].
\end{split}
\end{equation}

For the $\WADD$ we have by \eqref{eq:N_Lambda_same}, 
\eqref{eq:sum_lambdas_equal} and \eqref{eq:GCuSum_sumoflambda} for each $\theta \in \Theta$,
\begin{equation}\label{eq:GDECuSum_WADDProofStmnt}
\begin{split}
\Expect^\theta_1[\taugd(w)]
&=\Expect^\theta_1\left[\sum_{n=1}^{\bar{N}(w)}\bar{\lambda}_k(w)\right] + \Expect^\theta_1\left[\sum_{n=1}^{\bar{N}(w)-1}t_k(w)\right]\\
&=\Expect^\theta_1\left[\sum_{n=1}^{N(w)} \lambda_k(w)\right] + \Expect^\theta_1\left[\sum_{n=1}^{\bar{N}(w)-1}t_k(w)\right]\\
&=\Expect^\theta_1\left[\taugc \right] + \Expect^\theta_1\left[\sum_{n=1}^{\bar{N}(w)-1}t_k(w)\right]\\
&\leq\Expect^\theta_1\left[\taugc \right] + \Expect^\theta_1\left[\bar{N}(w)-1\right] \lceil h/\mu \rceil \\
&=\Expect^\theta_1\left[\taugc \right] + \Expect^\theta_1\left[N(w)-1\right] \lceil h/\mu \rceil \\
&\leq \Expect^\theta_1[\taugc] + \frac{1}{q_\theta} \lceil h/\mu \rceil.
\end{split}
\end{equation}
In \eqref{eq:GDECuSum_WADDProofStmnt} we have also used the fact that
\[
t_k(w) \leq \lceil h/\mu \rceil, \; \forall w\in [0,A), \forall k,
\]
and the upper bound obtained on $\Expect^\theta_1[N(w)]$ from \eqref{eq:ENx_UB_11}.
Also note that the right-hand side is not a function of $w$, but does depend on the assumption that $w\in [0,A)$.

We now obtain an upper bound on $\Expect^\theta_\gamma[(\taugd-\gamma)^+|\mathcal{I}_{\gamma-1}]$.
If $\mathcal{I}_{\gamma-1}=i_{\gamma-1}$ is such that $W_{\gamma-1}(\theta^*)=w \in [0,A)$, then
\[
\Expect^\theta_\gamma[(\taugd-\gamma)^+|\mathcal{I}_{\gamma-1}=i_{\gamma-1}] \leq \Expect^\theta_1[\taugd(w)].
\]
This is because for $n \geq \gamma$
\begin{equation}\label{eq:GLRT_related_stmnt}
\max_{1\leq k \leq n} \;\sup_{\theta \in \Theta} \; \sum_{i=k}^n \log \frac{f_\theta(X_i^{(S_i)})}{f_0(X_i^{(S_i)})}
\geq \max_{\gamma \leq k \leq n} \;\sup_{\theta \in \Theta} \; \sum_{i=k}^n \log \frac{f_\theta(X_i^{(S_i)})}{f_0(X_i^{(S_i)})}.
\end{equation}
Thus, if $\mathcal{I}_{\gamma-1}=i_{\gamma-1}$ is such that $W_{\gamma-1}(\theta^*)=w \in [0,A)$, then using
\eqref{eq:GDECuSum_WADDProofStmnt} we have
\begin{equation}\label{eq:GDECuSum_WADDProofStmnt_1}
\begin{split}
\Expect^\theta_\gamma[(\taugd-\gamma)^+|\mathcal{I}_{\gamma-1}=i_{\gamma-1}] \leq \Expect^\theta_1[\taugd(w)]
\leq \Expect^\theta_1[\taugc] + \frac{1}{q_\theta}  \lceil h/\mu \rceil.
\end{split}
\end{equation}

On the other hand,
if $\mathcal{I}_{\gamma-1}=i_{\gamma-1}$ is such that $W_{\gamma-1}(\theta^*)=w <0$, then
the time to cross $A$ for the GDECuSum
statistic will be equal to the time taken for the statistic $W_{\gamma-1}(\theta^*)$ to cross $0$ from below, plus
a time bounded by $\Expect_1[\taugd(0)]$, where again we have used \eqref{eq:GLRT_related_stmnt}.
 Thus, we can write,
\begin{equation}\label{eq:GDECuSum_WADDProofStmnt_2}
\begin{split}
\Expect^\theta_\gamma[(\taugd-\gamma)^+|\mathcal{I}_{\gamma-1}=i_{\gamma-1}]
&\leq \lceil h/\mu \rceil + \Expect^\theta_1[\taugd(0)] \\
&\leq \Expect^\theta_1[\taugc] + (\frac{1}{q_\theta}+1)  \lceil h/\mu \rceil.
\end{split}
\end{equation}
Thus, we can write due to \eqref{eq:CADDWADDWorstGamma_MCuSum}
\begin{equation}\label{eq:GDECuSum_WADDProofStmnt_3}
\begin{split}
\Expect^\theta_\gamma[(\taugd-\gamma)^+|\mathcal{I}_{\gamma-1}]
&\leq \Expect^\theta_1[\taugc] + (\frac{1}{q_\theta}+1)  \lceil h/\mu \rceil\\
&= \WADD^\theta(\taugc) + (\frac{1}{q_\theta}+1)  \lceil h/\mu \rceil + 1.
\end{split}
\end{equation}
Note that the right-hand side is no more a function of the conditioning $\mathcal{I}_{\gamma-1}$.
The proof is complete if we define
\[
K_{\text{GD}} = (\frac{1}{q_\theta}+1)  \lceil h/\mu \rceil + 1,
\]
and take the essential supremum on the left-hand side.
\end{IEEEproof}

\section*{Acknowledgements}
We thank Dr. Sirin Nitinawarat for discussions that led to the work described in Section~\ref{sec:ExistOfLeastFavDist}.

%\pagebreak
%\nocite{*}
\bibliographystyle{ieeetr}

%\bibliographystyle{elsarticle-harv}

%\bibliographystyle{elsarticle-harv}
%\bibliography{QCD_verVV}
\bibliography{QCD_verSubmitted}
\end{document}